\input amstex
\documentstyle{amsppt}
\document
\topmatter
\title
Solvable points on genus one curves over local fields
\endtitle
\title
Solvable points on genus one curves over local fields
\endtitle
\author
Ambrus P\'al
\endauthor
\date
February 8, 2012.
\enddate
\address
Department of Mathematics, 180 Queen's Gate, Imperial College,
London SW7 2AZ, United Kingdom\endaddress
\email a.pal\@imperial.ac.uk\endemail
\abstract Let $F$ be a field complete with respect to a discrete valuation whose residue field is perfect of characteristic $p>0$. We prove that every smooth, projective, geometrically irreducible curve of genus one defined over $F$ with a non-zero divisor of degree a power of $p$ has a solvable point over $F$. We also show that there is a field $F$ complete with respect to a discrete valuation whose residue field is perfect and there is a finite Galois extension $K|F$ such that there is no solvable extension $L|F$ such that the extension $KL|K$ is unramified, where $KL$ is the composite of $K$ and $L$. As an application we deduce that that there is a field $F$ as above and there is a smooth, projective, geometrically irreducible curve over $F$ which does not acquire semi-stable reduction over any solvable extension of  $F$.
\endabstract
\footnote" "{\it 2000 Mathematics Subject Classification. \rm 14G05, 11G07.}
\endtopmatter

\heading 1. Introduction
\endheading

We say that a finite Galois field extension $K|F$ is solvable if its Galois group is solvable. A field $F$ is solvably closed if it has no proper finite solvable Galois extensions. For every field $F$ let $F^{per}$ denote its perfection and let $F^{sol}$ denote the maximal Galois extension of $F^{per}$ with pro-solvable Galois group over $F^{per}$. We will call $F^{sol}$ the solvable closure of $F$. As its name indicates $F^{sol}$ is solvably closed. Let $X$ be a quasi-projective variety over a field $K$. By a solvable point of $X$ over $F$ we mean a $F^{solv}$-valued point of $X$. We say that $X$ has solvable points if $X(F^{solv})$ is non-empty. The main motivation for writing this paper is to give a correct proof to the following
\proclaim{Theorem 1.1} Let $F$ be a field complete with respect to a discrete valuation whose residue field is perfect of characteristic $p>0$. Then every smooth, projective, geometrically irreducible curve of genus one defined over $F$ with a non-zero divisor of degree a power of $p$ has a solvable point over $F$. 
\endproclaim
The motivation for Theorem 1.1 is my conjecture saying that every genus one curve defined over any field has a solvable point. The theorem above is the only one other known non-trivial case of this conjecture besides the theorem of \c{C}iperiani-Wiles (see [1]). Theorem 1.1 was announced in my paper [14] which contained two proofs. The first, more complete proof rests crucially on Lemma 6.4 on page 632, which is false. In fact in this paper we will show the following:
\proclaim{Theorem 1.2} There is a field $F$ of characteristic zero complete with respect to a discrete valuation whose residue field is perfect and there is a finite Galois extension $K|F$ such that there is no solvable extension $L|F$ such that the extension $KL|K$ is unramified, where $KL$ is the composite of $K$ and $L$.
\endproclaim
The paper [14] also contains the sketch of a second strategy. This argument can be completed, although one needs significantly more sophisticated tools, for example the Oort-Tate classification theorem, than what is mentioned in my original article. In the course of the proof of Theorem 1.1 we will also prove a result (Theorem 3.6) which is interesting on its own, and which gives a rather explicit description of the $p$-torsion cohomology of certain elliptic curves defined over local fields whose residue field is perfect of characteristic $p$. Moreover for the above-mentioned reasons this paper will also contain some errata. 

The proof of Theorem 1.2 is also quite intricate, for example it uses some non-trivial facts about Drinfeld modular curves. However it is mainly interesting because it has the following
\proclaim{Corollary 1.3} There is a field $F$ of characteristic zero complete with respect to a discrete valuation whose residue field is perfect and there is a smooth, projective, geometrically irreducible curve over $F$ which does not acquire semi-stable reduction over any solvable extension of $F$.
\endproclaim
This result is significant because it opens up the possibility to go beyond the results of my paper [14] in its quest to construct curves of a given genus without solvable points. In a joint paper [8] Gyula K\'arolyi and I determined the set of natural numbers for which the original method of [14] can work which uses curves with semi-stable reduction. This set does not contain the number $23$, for example, but it looks unlikely that every genus $23$ curve has a solvable point. However it might be possible to construct counter-examples with this genus by using a curve which will not acquire semi-stable reduction over any solvable extension. 

The contents of this paper are the following. In the next section we recall the Oort-Tate classification [13] of finite flat group schemes of order $p$ over certain complete local rings and a theorem of Roberts (see [15]) which computes the cohomology of these groups schemes in terms of this classification (in the {\it fppf} cohomology). In the third section we use these results to give an explicit description of the cohomology of the $p$-torsion of certain elliptic curves defined over such base rings. In the fourth section first we will give a short proof of the fact that an elliptic curve defined over a local field with perfect residue field will attain semi-stable extension after a finite solvable extension then we use our previous results to prove Theorem 1.1. We prove Theorem 1.2 in the fifth section. In the last section we first give an overview of Gerritzen's uniformisation theory of abelian varieties by non-archimedean tori, then we derive Corollary 1.3.
\definition{Errata for [14]} Currently we only know Theorem 6.3 of [14] on page 632 when the residue field is perfect. Consequently Theorem 6.6 of [14] on page 632 is known under this assumption only. Lemma 6.4 of [14] on page 632 is false. Here we need to assume that the extension $L|F$ in the claim is tamely ramified, otherwise the conclusion does not hold in general (see Theorem 1.2). With this additional assumption the claim is true (see Lemma 4.6). Moreover the remark on the vanishing of the cohomology group $H^1(K,E_n)$ on page 637 of [14] is also false (see Theorem 3.6).
\enddefinition
\definition{Acknowledgement 1.4} The author was partially supported by the EPSRC grant P19164. 
\enddefinition

\heading 2. The cohomology of Oort-Tate group schemes
\endheading

\definition{Notation 2.1} We say that a field $F$ is local if it is complete with respect to a discrete valuation. For the sake of simplicity we will assume that every local field in this paper has perfect residue field. For every local field $F$ let $\Cal O_F$ and $\bold k_F$ denote its valuation ring and its residue field, respectively. Let $F$ be a local field of characteristic zero whose residue field $\bold k_F$ has characteristic $p>0$. Let $\text{ord}_F$ denote the valuation of $F$ normalised such that $\text{ord}_F(\pi)=1$ for every uniformizer $\pi\in\Cal O_F$. Let $e=\text{ord}_F(p)$ denote the absolute ramification index of $F$. For every finite flat group scheme $G$ over any base let $G^D$ denote its Cartier dual. For every discrete valuation ring $R$ and for every finite flat group scheme $G$ over Spec$(R)$ let Disc$(G)$ denote the discriminant ideal of $G$, which is by definition the discriminant ideal of the finite flat $R$-algebra $\Gamma(G,\Cal O_G)$ over $R$. Such an ideal is uniquely determined by its order, that is, the valuation of any of its generators. For every $G$ as above which is annihilated by $p$ and for every $s\in\Bbb F_p^*$ let $[s]:\Gamma(G,\Cal O_G)\rightarrow\Gamma(G,\Cal O_G)$ denote the ring homomorphism induced by the multiplication by $s$ map on $G$. Let $\omega:\Bbb F_p^*\rightarrow\Cal O_F^*$ denote the Teichm\"uller character. 
\enddefinition
\proclaim{Theorem (Oort-Tate) 2.2} The following holds:
\roster
\item"$(i)$" For every $a\in\Cal O_F$ which divides $p$ there is a unique finite flat group scheme $G_a$ over $\text{\rm Spec}(\Cal O_F)$ such that as a scheme
$$G_a=\text{\rm Spec}(\Cal O_F[x]/(x^p-ax))$$
and for every $s\in\Bbb F_p^*$ we have $[s]x=\omega(s)x$,
\item"$(ii)$" every finite flat commutative group scheme $G$ over $\text{\rm Spec}(\Cal O_F)$ of rank $p$ is isomorphic to $G_a$ for some $a\in\Cal O_F$ as above,
\item"$(iii)$" there is a non-trivial homomorphism $\phi:G_a\rightarrow G_b$ of group schemes if and only if there is an element $u\in\Cal O_F$ such that $au^{p-1}=b$. This homomorphism is an isomorphism if and only if $u\in\Cal O_F^*$,
\item"$(iv)$" we have $\text{\rm Disc}(G_a)=(a^p)$ and $G_a^D=G_{wp/a}$ for an element $w\in\Cal O^*_F$ which does not depend on $G$.
\endroster
\endproclaim
\definition{Proof} Claims $(i)$, $(ii)$ and the first half of claim $(iii)$ are special cases of Theorem 4.4.1 of [17] on page 149. Assume now that $\phi:G_a\rightarrow G_b$ is an isomorphism of group schemes. By the above there is an element $u\in\Cal O_F$ such that $au^{p-1}=b$. Since the inverse $\phi^{-1}:G_b\rightarrow G_a$ of $\phi$ is also a non-trivial homomorphism there is an element $v\in\Cal O_F$ such that $bv^{p-1}=a$. Clearly $u^{-1}=v$ and hence $u\in\Cal O_F^*$. Assume now there is an element $u\in\Cal O_F^*$ such that $au^{p-1}=b$. Write
$$G_a=\text{\rm Spec}(\Cal O_F[x]/(x^p-ax)),\quad 
G_b=\text{\rm Spec}(\Cal O_F[y]/(y^p-by))$$
such that for every $s\in\Bbb F_p^*$ we have $[s]x=\omega(s)x$ and $[s]y=\omega(s)y$. By claim $(iii)$ of Theorem 4.4.1 of [17] on page 149 there is a homomorphism $\phi:G_a\rightarrow G_b$ of group schemes such that
$\phi^*(y)=ux$. This map is an isomorphism of schemes; its inverse is the unique morphism $\psi:G_b\rightarrow G_a$ of schemes such that $\psi^*(x)=u^{-1}y$. The map $\psi$ is a homomorphism of group schemes by claim $(iii)$ of Theorem 4.4.1 of [17], so the second half of claim $(iii)$ is true. The first half of claim $(iv)$ is an immediate consequence of claim $(i)$ while the second half of claim $(iv)$ is just the remark at the top of page 15 of [13].\qed 
\enddefinition
For every finite flat commutative group scheme $G$ over $\text{\rm Spec}(\Cal O_F)$ let $G_F$ denote its generic fibre, that is, the base change of $G$ to Spec$(F)$. We say that a $G$ as above splits generically when $G_F$ is a constant group scheme. In the next claim we assume that $F$ contains the $p$-th roots of unity. In this case $m=e/(p-1)$ is an integer.
\proclaim{Corollary 2.3} The following holds:
\roster
\item"$(i)$" the group scheme $G_a$ splits generically if and only if $a=u^{p-1}$ for some $u\in\Cal O_F$,
\item"$(ii)$" two finite flat groups schemes of rank $p$ over $\text{\rm Spec}(\Cal O_F)$ which split generically are isomorphic if and only if their discriminants are the same, 
\item"$(iii)$" the possible orders of the discriminants of such groups schemes are the integers $p(p-1)i$ where $0\leq i\leq m$.
\endroster
\endproclaim
\definition{Proof} This is the Corollary on page 689 of [15]. Since no proof can be found there, we include one here for the reader's convenience. Let $G$ denote the constant group scheme of order $p$ over Spec$(\Cal O_F)$. For every $\alpha\in\Bbb F_p=G(\text{Spec}(\Cal O_F))$ let $e_{\alpha}\in\Gamma(G,\Cal O_G)$ be the idempotent corresponding to the connected component of $G$ which is the image of the section $\alpha:\text{Spec}(\Cal O_F)\rightarrow G$. Let $y$ be the element:
$$y=\sum_{\alpha\in\Bbb F_p^*}\omega(\alpha)^{-1}e_{\alpha}\in\Gamma(G,\Cal O_G).$$
Then for every $\beta\in\Bbb F_p^*$ we have:
$$[\beta]y=\sum_{\alpha\in\Bbb F_p^*}\omega(\alpha)^{-1}e_{\alpha\beta}=
\omega(\beta)y.$$
Clearly $\Gamma(G,\Cal O_S)=\Cal O_F[y]$ and since
$$y^p=\sum_{\alpha\in\Bbb F_p^*}\omega(\alpha)^{-p}e_{\alpha}^p=
\sum_{\alpha\in\Bbb F_p^*}\omega(\alpha)^{-1}e_{\alpha}=y,$$
we get that $G_1$ is the constant group scheme of order $p$ over $\text{\rm Spec}(\Cal O_F)$. Let $a\in\Cal O_F$ be such that $a$ divides $p$ and $G_a$ splits generically. Then there is an isomorphism $\phi_0:(G_1)_F\rightarrow(G_a)_F$ of group schemes. Since as a scheme $G_1$ is the coproduct of copies of Spec$(\Cal O_F)$ there is a unique prolongation $\phi:G_1\rightarrow G_a$ of $\phi_0$ by the valuative criterion of properness. The map $\phi$ is a homomorphism of group schemes because it is such a map on a dense open subscheme. Hence $a=u^{p-1}$ for some $u\in\Cal O_F$ by claim $(iii)$ of Theorem 2.2. On the other hand if $u\in\Cal O_F$ is such that $a=u^{p-1}$ divides $p$ then there is a non-trivial homomorphism $\phi:G_1\rightarrow G_a$ of group schemes by claim $(iii)$ of Theorem 2.2. Such a homomorphism induces an isomorphism between $(G_1)_F$ and $(G_a)_F$ and hence $G_a$ spits generically. So claim $(i)$ is true. 
Let again $u\in\Cal O_F$ be such that $a=u^{p-1}$ divides $p$. By claim $(iv)$ of Theorem 2.2 the order of the discriminant of $G_a$ is $p\cdot\text{ord}_F(a)=p(p-1)\text{ord}_F(u)$. Since the value of $\text{ord}_F(u)$ can be any integer between $0$ and $\text{ord}_F(p)/(p-1)=m$, claim $(iii)$ is true. On the other hand if $v\in\Cal O_F$ is such that $b=v^{p-1}$ divides $p$ and the orders of the discriminants of the groups schemes $G_a$ and $G_b$ are the same then $u/v\in\Cal O_F^*$ by the above. Therefore $G_a$ and $G_b$ are isomorphic by claim $(iii)$ of Theorem 2.2, so claim $(ii)$ holds, too.\qed 
\enddefinition
\definition{Definition 2.4} In this paper by the cohomology of a group scheme $G$ over a base scheme $X$ we mean the cohomology of the sheaf it represents for the {\it fppf}-topology on $X$. Note that for a smooth quasi-projective group scheme this is the same as its cohomology with respect to the \'etale topology by Theorem 2.9 of [11] on page 114. For every $X$ of the form Spec$(A)$ where $A$ is a commutative ring we let $H^i(A,G)$ denote $H^i(\text{Spec}(A),G)$ for the sake of simple notation. Let $F$ be as above and assume that $F$ contains the $p$-th roots of unity. Let $G$ be a finite flat group scheme of rank $p$ over $\text{\rm Spec}(\Cal O_F)$ which splits generically. Fix an isomorphism $a:G_F\rightarrow\mu_p$ and let $\alpha$ denote the composition of the maps:
$$\CD H^1(\Cal O_F,G)@>>>H^1(F,G_F)@>>>H^1(F,\mu_p)\endCD$$
where the first map is furnished by base change and the second map is the isomorphism induced by $a$. Let
$$\delta:F^*/(F^*)^p\longrightarrow H^1(F,\mu_p)\tag2.4.1$$
be the coboundary map furnished by the Kummer short exact sequence:
$$\CD 0@>>>\mu_p@>>>\Bbb G_m@>{x\mapsto x^p}>>\Bbb G_m@>>>0,\endCD$$
where $\Bbb G_m$ denote the multiplicative group scheme. For the sake of simple notation let $U_F=\Cal O_F^*$ and for every $i\in\Bbb N$ let
$$U_F^{(i)}=\{u\in U_F|\text{ord}_F(1-u)\geq i\}.$$
\enddefinition
\proclaim{Theorem (Roberts) 2.5} Let $G$ be a finite flat group scheme of rank $p$ over $\text{\rm Spec}(\Cal O_F)$ which splits generically. Then there is a commutative diagram:
$$\CD H^1(\Cal O_F,G)@>{\alpha}>>H^1(F,\mu_p)\\
@A{\beta}AA@AA{\delta}A\\
U_F^{(i)}U_F^p/U_F^p@>{\iota}>> F^*/(F^*)^p\endCD$$
where $\beta$ is an isomorphism, the homomorphism $\iota$ is the inclusion map and
$$i={pe-\text{\rm ord}(\text{\rm Disc}(G))\over p-1}.$$
\endproclaim
\definition{Proof} This is Theorem 1 of [15] on page 694. Note however that the claim there contains a typographical error; compare it with the announcement of this result on page 228 of [10]. Also note that although the result is formulated for fields with finite residue fields only, the proof works in general. In fact the argument in [15] which shows that $\alpha$ is injective (Proposition 3 on page 692) works for any discrete valuation ring, the proof of the fact that the image of the map $\alpha$ lies in the image of $\delta\circ\iota$ (see page 694) only uses standard facts about discriminants, while the second proof of the surjectivity of $\alpha$ onto the image of $\delta\circ\iota$ (see pages 698--700) is an explicit construction.\qed 
\enddefinition
\heading 3. The cohomology of the $p$-torsion of admissible elliptic curves
\endheading
\definition{Definition 3.1} For every commutative group scheme $A$ and for every positive integer $m$ let $A[m]$ denote the $m$-torsion sub-group scheme of $A$. For every field $K$ let $\overline K$ denote its separable closure. We continue to denote by $F$ a local field of characteristic zero with a residue field of characteristic $p>0$. For every such $F$ let $\Gamma_F$ and $I_F$ denote the absolute Galois group of $F$ and the inertia subgroup of $\Gamma_F$, respectively. We say that an elliptic curve $E$ over $F$ is {\it admissible} if the action of $\Gamma_F$ on $E[p]$ is unramified and does not factor through a solvable quotient. Let $E$ be an admissible elliptic curve over $F$ and let $F(E)$ be the field of definition of the $\overline F$-valued points of $E[p]$. The field $F(E)$ contains the $p$-th roots of unity. However the extension $F(E)|F$ is unramified while the extension of $F$ which we get by adjoining the $p$-th roots of unity is totally ramified, therefore $F$ must contain the $p$-th roots of unity in this case.
\enddefinition
\proclaim{Proposition 3.2} Assume that the admissible elliptic curve $E$ over $F$ has good reduction and let $\Cal E$ denote its N\'eron model over $\text{\rm Spec}(\Cal O_F)$. Then the base change of the finite flat group scheme $\Cal E[p]$ to $\text{\rm Spec}(\Cal O_{F(E)})$ is isomorphic to $G^{\oplus2}$ where $G$ is the unique finite flat group scheme of rank $p$ over $\text{\rm Spec}(\Cal O_{F(E)})$ which splits generically and the order of its discriminant is $p(p-1)m/2$.
\endproclaim
\definition{Proof} Let $G_0$ be a closed subgroup scheme of rank $p$ of the base change $H_0$ of $\Cal E[p]$ to $\text{\rm Spec}(F(E))$ and let $G$ be the scheme-theoretical closure of $G_0$ in the base change $H$ of $\Cal E[p]$ to $\text{\rm Spec}(\Cal O_{F(E)})$. Then $G$ is the prolongation of $G_0$, that is, a closed subgroup scheme of rank $p$ of $H$ whose generic fibre is $G_0$. By assumption there is an element $\sigma\in\Gamma(F(E)|F)$ such that the pull-back $\sigma^*(G_0)$ is a subgroup scheme of $H_0$ complementary to $G_0$. Because the extension $F(E)|F$ is unramified $\sigma$ is the extension of an automorphism of $\Cal O_{F(E)}$ which fixes $\Cal O_F$. The pull-back of $H$ with respect to this automorphism is isomorphic to $H$, and hence the prolongation of $\sigma^*(G_0)$ in $H$ is isomorphic to $G$. The direct sum of these two copies of $G$ has the same rank as $H$ therefore these group schemes are equal. The group scheme $H$, being the $p$-torsion of an elliptic curve, is self-dual, so $G^{\oplus2}$ is isomorphic to $(G^D)^{\oplus2}$. Suitable projection maps furnish non-trivial homomorphisms $G\rightarrow G^D$ and $G^D\rightarrow G$ of group schemes. Therefore $G$ and $G^D$ are isomorphic by claim $(iii)$ of Theorem 2.2. Choose an $a\in\Cal O_F$ dividing $p$ such that $G$ is isomorphic to $G_a$. Then $G^D$ is isomorphic to $G_{wp/a}$ for a unit $w\in\Cal O_F^*$ by claim $(iv)$ of Theorem 2.2. Therefore the orders of the discriminants of $G_a$ and $G_{wp/a}$ are the same, so 
$$\text{ord}_F(a)=\text{ord}_F(p)-\text{ord}_F(a)=(p-1)m-\text{ord}_F(a),$$
and hence the order of Disc$(G)$ is $p(p-1)m/2$ by claim $(iv)$ of Theorem 2.2.\qed 
\enddefinition
\definition{Definition 3.3} For every finite extension $L|K$ of fields let $\Gamma(L|K)$ denote its automorphism group. Fix an isomorphism $\Bbb F_p\rightarrow\mu_p$ of group schemes over Spec$(F)$. The cup product induces a $\Gamma(F(E)|F)$-equivariant map:
$$\cup:H^1(F,\Bbb F_p)\otimes H^0(F(E),E[p])\longrightarrow H^1(F(E),E[p]).\tag3.3.1$$
By slight abuse of notation let $E[p]$ also denote the $\Gamma(L|K)$-module $H^0(F(E),E[p])$. Let 
$$\delta:F(E)^*/(F(E)^*)^p\otimes E[p]\longrightarrow H^1(F(E),E[p])$$
denote the composition of the tensor product of the coboundary map in (2.4.1) and the identity map of $E[p]$ with the map in (3.3.1). Clearly this map is an isomorphism of $\Gamma(F(E)|F)$-modules. We will equip the tensor product of $\Gamma(F(E)|F)$-modules with the usual $\Gamma(F(E)|F)$-action.
\enddefinition
\proclaim{Proposition 3.4} Let $E$ be an admissible elliptic curve over $F$ with good reduction and let $\Cal E$ denote its N\'eron model over $\text{\rm Spec}(\Cal O_F)$. Then there is a commutative diagram:
$$\CD H^1(\Cal O_F,\Cal E[p])@>{\alpha}>>H^1(F,E[p])\\
@A{\beta}AA@AA{\delta}A\\
\left(U_{F(E)}^{(i)}U_{F(E)}^p/U_{F(E)}^p\otimes E[p]\right)^{\Gamma(F(E)|F)}@>{\iota}>>
\left(F(E)^*/(F(E)^*)^p\otimes E[p]\right)^{\Gamma(F(E)|F)}\endCD$$
where $\alpha$ is furnished by base change, the map $\beta$ is an isomorphism, the homomorphism $\iota$ is the inclusion map, and $i=pm/2$.
\endproclaim
\definition{Proof} Let $G$ be the finite flat group scheme over Spec$(\Cal O_{F(E)})$ introduced in Proposition 3.2. Clearly $H^1(\text{Spec}(\Cal O_{F(E)}),\Cal E[p])=H^1(\text{Spec}(\Cal O_{F(E)}),G)^{\oplus2}$ so by Theorem 2.5 there is a commutative diagram:
$$\CD H^1(\Cal O_{F(E)},\Cal E[p])@>>>H^1(F(E),E[p])\\
@A{\beta}AA@AA{\delta}A\\
U_{F(E)}^{(i)}U_{F(E)}^p/U_{F(E)}^p\otimes E[p]
@>{\iota}>>
F(E)^*/(F(E)^*)^p\otimes E[p]\endCD$$
where $\beta$ is an isomorphism, the homomorphism $\iota$ is the inclusion map, and $i=pm/2$. Note that $\Gamma(F(E)|F)$ acts on all groups above by functoriality and every map in the diagram is equivariant with respect to these actions. Therefore we get another commutative diagram:
$$\CD H^1(\Cal O_{F(E)},\Cal E[p])^{\Gamma(F(E)|F)}
@>>>H^1(F(E),E[p])^{\Gamma(F(E)|F)}\\
@A{\beta}AA@AA{\delta}A\\
\left(U_{F(E)}^{(i)}U_{F(E)}^p/U_{F(E)}^p\otimes E[p]\right)^{\Gamma(F(E)|F)}@>{\iota}>>
\left(F(E)^*/(F(E)^*)^p\otimes E[p]\right)^{\Gamma(F(E)|F)}\endCD$$
by taking $\Gamma(F(E)|F)$-invariants. By Lemma 6.8 of [14] on page 635 the cohomology groups $H^r(\Gamma(F(E)|F),E[p])$ vanish for every $r$. Therefore every $E_{r,0}$ term of the Hochschild-Serre spectral sequence $$H^r(\Gamma(F(E)|F),H^s(F(E),E[p]))\Rightarrow H^{r+s}(F,E[p])$$
vanishes, and hence the restriction map $H^1(F,E[p])\rightarrow
H^1(F(E),E[p])^{\Gamma(F(E)|F)}$ is an isomorphism. By part $(a)$ of Remark 2.21 of [11] on page 105 there is a spectral sequence:
$$H^r(\Gamma(F(E)|F),H^s(\Cal O_{F(E)},\Cal E[p]))\Rightarrow H^{r+s}(\Cal O_F,E[p]).$$
By the valuative criterion of properness $H^0(F(E),E[p])=H^0(\Cal O_{F(E)},\Cal E[p])$. Therefore the same argument as above implies that the restriction map $H^1(\Cal O_F,\Cal E[p])\rightarrow H^1(\Cal O_{F(E)},\Cal E[p])^{\Gamma(F(E)|F)}$ is an isomorphism, too. Now the claim follows from the naturality of the restriction maps.\qed 
\enddefinition
\definition{Definition 3.5} Let $F$ be as above and let $E$ be an admissible elliptic curve over $F$. We say that a cohomology class $c\in H^1(F,E[p])$ is small if its image under the restriction map:
$$j_{E/F}:H^1(F,E[p])\longrightarrow H^1(F(E),E[p])^{\Gamma(F(E)|F)}$$
lies in the image of the map:
$$\delta\circ\iota:U_{F(E)}^{(mp/2)}U_{F(E)}^p/U_{F(E)}^p\otimes E[p]\longrightarrow H^1(F(E),E[p])$$
where we continue to let $m$ denote $e/(p-1)$ where $e$ is the absolute ramification index of $F$. Let $\iota_{E/F}:H^1(F,E[p])\rightarrow H^1(F,E)$ be the map induced by the inclusion $E[p]\subset E$.
\enddefinition
\proclaim{Theorem 3.6} Let $E$ be an admissible elliptic curve over $F$ with good reduction and let $c\in H^1(F,E[p])$ be a cohomology class. Then $\iota_{E/F}(c)\in H^1(F,E)[p]$ is zero if and only if $c$ is small.
\endproclaim
\definition{Proof} Let $\Cal E$ be the N\'eron model of $E$ over Spec$(\Cal O_F)$. Then we have the following commutative diagram:
$$\CD\Cal E(\Cal O_F)@>>>H^1(\Cal O_F,\Cal E[p])@>>>H^1(\Cal O_F,
\Cal E)[p]\\
@VVV@VVV@VVV\\
E(F)@>>>H^1(F,E[p])@>\iota_{E/F}>>H^1(F,E)[p]\endCD$$
where the upper and lower rows are parts of the cohomological long exact sequences furnished by the Kummer exact sequences of $\Cal E$ and $E$, respectively. By the N\'eron extension property the left vertical arrow is an isomorphism. Therefore $c$ is small if $\iota_{E/F}(c)$ is zero, by Proposition 3.4. In order to show that the converse also holds it will be sufficient to show the following
\enddefinition
\proclaim{Lemma 3.7} We have $H^1(\Cal O_F,\Cal E)[p]=0$.
\endproclaim
\definition{Proof} Let $d\in H^1(\Cal O_F,\Cal E)$ be a $p$-torsion cohomology class. Let $T\rightarrow\text{Spec}(\Cal O_F)$ be an $\Cal E$-torsor which represents $d$. Since every torsor over $\Cal E$ is smooth over $\text{Spec}(\Cal O_F)$ by Proposition 4.2 of [11] on page 120, it will be sufficient to show that $T$ has a $\bold k_F$-valued point by Hensel's lemma. In order to do so it will be sufficient to prove that $H^1(\bold k_F,\Cal E_0)[p]=0$ where $\Cal E_0$ denotes the the special fibre of $\Cal E$. By Proposition 3.2 the base change of the finite flat group scheme $\Cal E_0[p]$ to $\text{\rm Spec}(\overline{\bold k}_F)$ is isomorphic to $G_0^{\oplus2}$ where $G_0$ is the base change of the group scheme $G$ of Proposition 3.2 to Spec$(\overline{\bold k}_F)$. This is only possible when $\Cal E_0$ is a supersingular elliptic curve. Therefore the $p$-power map $\Cal E_0(\overline{\bold k}_F)\rightarrow\Cal E_0(\overline{\bold k}_F)$ given by the rule $x\mapsto x^p$ is an isomorphism. The induced map $H^1(\bold k_F,\Cal E_0)\rightarrow H^1(\bold k_F,\Cal E_0)$ is multiplication by $p$. It is also an isomorphism, so its kernel is trivial.\qed
\enddefinition
\heading 4. Solvable points on genus one curves over local fields
\endheading
\definition{Notation 4.1} For every variety $V$ defined over a field $F$ and for every extension $K$ of $F$ let $V_K$ denote the base change of $V$ to Spec$(K)$. For every abelian variety $A$ defined over a field $F$ and for every prime number $l$ different from the characteristic of $F$ let $T_l(A)$ be the $l$-adic Tate module of $A$ and let $V_l(A)$ denote the Galois representation $T_l(A)\otimes_{\Bbb Z_l}\Bbb Q_l$ over $F$.
\enddefinition
Recall that we say that an abelian variety $A$ defined over a local field $F$ has semi-stable reduction if the fibre over Spec$(\bold k_F)$ of its N\'eron model over Spec$(\Cal O_F)$ is a semi-abelian variety. We will need the following two fundamental results.
\proclaim{Theorem 4.2} Let $F$ be a local field and let $A$ be an abelian variety defined over $F$. Then there is a finite Galois extension $K$ of $F$ such that $A_K$ has semi-stable reduction over $K$.
\endproclaim
\definition{Proof} This is Th\'eor\`eme 3.6 of [6] on page 351.\qed
\enddefinition
\proclaim{Theorem 4.3} Let $F$ be a local field and let $A$ be an abelian variety defined over $F$. Let $l$ be a prime number different from the characteristic of $\bold k_F$. Then the following conditions are equivalent:
\roster
\item"$(i)$" the action of $I_F$ on $V_l(A)$ is unipotent,
\item"$(ii)$" the abelian variety $A$ has semi-stable reduction over $F$.
\endroster
\endproclaim
\definition{Proof} This follows from Proposition 3.5 of [6] on page 350
and Corollaire 3.8 of [6] on page 353.\qed
\enddefinition
\proclaim{Proposition 4.4} Let $F$ be a local field and let $E$ be an elliptic curve over $F$. Then there is a finite solvable Galois extension $K$ of $F$ such that $E_K$ has semi-stable reduction over $K$.
\endproclaim
\definition{Proof} Let $l=2$ if the characteristic of $\bold k_F\neq2$ and let $l=3$, otherwise. In either case the automorphism group of the group $E[l](\overline F)$ is solvable hence we may assume without the loss of generality, by taking a finite solvable extension of $F$ if it is necessary, that the action of $\Gamma_F$ on $E[l](\overline F)$ is trivial. In this case the image of $\Gamma_F$ in the automorphism group of the $\Bbb Z_l$-module $T_l(E)$ is a pro-$l$ group. By Theorem 4.2 there is a finite Galois extension $L$ of $F$ such that $E_L$ has semi-stable reduction over $L$. By the above the image of $\Gamma_L$ in the automorphism group of the $\Bbb Z_l$-module $T_l(E)$ is a normal subgroup of the image of $\Gamma_F$ whose index is a power of $l$. Hence there is a finite Galois extension $K|F$ such that $E_K$ has semi-stable reduction over $K$ and the group Gal$(K|F)$ is an $l$-group by Theorem 4.3. Because every $l$-group is solvable the claim is now clear.\qed
\enddefinition
\proclaim{Corollary 4.5} Let $F$ be a local field of characteristic zero whose residue field $\bold k_F$ has characteristic $p>0$. Let $E$ be an elliptic curve over $F$ and assume that the action of $\Gamma_F$ on $E[p]$ does not factor through a solvable quotient. Then there is a finite solvable Galois extension $K$ of $F$ such that $E_K$ has good supersingular reduction over $K$.
\endproclaim
\definition{Proof} We may assume, by taking a finite solvable extension of $F$ if it is necessary, that $E$ has semi-stable reduction, by Proposition 4.4. Because the stabiliser of a proper subgroup of $\Bbb F_p^2$ in $GL_2(\Bbb F_p)$ is solvable, the $\Bbb F_p$-linear action of $\Gamma_F$ on $E[p](\overline F)$ must be irreducible. Let $\Cal E$ be the N\'eron model of $E$ over Spec$(\Cal O_K)$ and let $\Cal E_0$ be the connected component of the identity of the fibre of $\Cal E$ over Spec$(\bold k_F)$. Because $E$ has semi-stable reduction by our assumption the group scheme $\Cal E_0$ is either isomorphic to the multiplicative group scheme $\Bbb G_m$ over a quadratic extension of $\bold k_F$ or $E$ has good reduction and $\Cal E_0$ is an elliptic curve. In the first case $E$ has Tate uniformisation over a quadratic extension of $F$. The image of the $p$-torsion of $\Bbb G_m$ with respect to this uniformisation is a proper $\Gamma_F$-submodule of $E[p](\overline F)$. So $\Cal E_0$ is an elliptic curve. By Hensel's lemma $E[p](\overline F)$ has a Galois submodule which is isomorphic to $\Cal E_0[p](\overline{\bold k}_F)$ as a group. Therefore $\Cal E_0$ is a supersingular elliptic curve since $\Cal E_0[p](\overline{\bold k}_F)$ would have order $p$ otherwise. The claim is now clear.\qed 
\enddefinition
\proclaim{Lemma 4.6} Let $F$ be a local field of characteristic $0$ and let $L|F$ be a finite tame Galois extension. Then there is a solvable extension $K|F$ such that the extension $KL|K$ is unramified, where $KL$ is the composite of $K$ and $L$.
\endproclaim
\definition{Proof} Let $T$ be the largest unramified subextension of $F$ in the field $L$. Then the Galois group Gal$(L|T)$ is cyclic of order $m$, where $m$ is relatively prime to the characteristic of $\bold k_F$. We may assume that $F$ contains the $m$-th roots of unity by adjoining them, if it is necessary. Then by Kummer theory the extension $L|T$ is the spitting field of a polynomial $x^m-u\pi^k$, where $k$ is an integer, $\pi$ is a uniformizer of $F$ and $u\in\Cal O_T^*$. Let $K$ be the splitting field of $x^m-\pi^k$ over $F$. Then the extension $LK|K$ is unramified since it is the composite of $T$ and the splitting field of $x^m-u$ over $K$.\qed
\enddefinition
\proclaim{Proposition 4.7} Let $F$ be a local field whose residue field $\bold k_F$ has characteristic $p>0$. Let $E$ be an elliptic curve over $F$ such that the action of $\Gamma_F$ on $E[p]$ does not factor through a solvable quotient. Then there is a finite solvable Galois extension $K$ of $F$ such that $E_K$ has good reduction over $K$ and it is admissible.
\endproclaim
\definition{Proof} Fix an isomorphism $E[p](\overline F)\cong\Bbb F_p^2$ of groups and let $G\leq GL_2(\Bbb F_p)$ be the image of $\Gamma_F$ under this identification. By Corollary 4.5 there is a finite solvable Galois extension $K$ of $F$ such that $E_K$ has good reduction over $K$. By the N\'eron--Ogg--Shafarevich condition $E_L$ will have good reduction over any finite extension $L$ of $F$ containing $K$, so we only need to show that there is a finite solvable Galois extension $L$ of $F$ such that $E_L$ is admissible. By taking an abelian extension, if it is necessary, we may assume that $G$ lies in $SL_2(\Bbb F_p)$. If $p$ does not divide the order of $G$ then the extension $F(E)|F$ is tame, and the claim holds by Lemma 4.6 in this case. Therefore we may assume that $p$ divides $|G|$. In this case any $p$-Sylow subgroup of $G$ is cyclic of order $p$, and it is also a $p$-Sylow of $SL_2(\Bbb F_p)$. By the Sylow theorem the number of $p$-Sylows of $G$ is congruent to $1$ modulo $p$, so $G$ either contains exactly one $p$-Sylow subgroup, or it contains all $p$-Sylows of $SL_2(\Bbb F_p)$, since the latter has $p+1$ copies of $p$-Sylow subgroups. In the former case $G$ is contained in a Borel subgroup of $SL_2(\Bbb F_p)$, the normaliser of some $p$-Sylow of $SL_2(\Bbb F_p)$, so it is solvable. In the latter case $G$ is a subgroup generated by the conjugate $p$-Sylows, and hence it is equal to $SL_2(\Bbb F_p)$ by Dickson's theorem (Theorem 8.4 of [5] on page 44). Let $H$ be the image of $I_F$ in $G$. We know that $p\geq 5$ since the group $GL_2(\Bbb F_p)$ is solvable otherwise. In this case $PSL_2(\Bbb F_p)$ is a simple group. By repeating the argument above which uses Dickson's theorem one can see that every normal subgroup of $SL_2(\Bbb F_p)$ either lies in the centre of $SL_2(\Bbb F_p)$ or it is equal to $SL_2(\Bbb F_p)$. Since $H$ is a normal subgroup we get that $H$ lies in the centre of $SL_2(\Bbb F_p)$ as $H$ is solvable and hence the second case is not possible. Therefore $p$ does not divide $|H|$ and hence the extension $F(E)|F$ is tame. So the claim holds by Lemma 4.6.\qed 
\enddefinition
\definition{Notation 4.8} For every finite extension $L$ of $F$ let $\text{ord}_L$ denote the valuation of $L$ normalised such that $\text{ord}_L(\pi)=1$ for every uniformizer $\pi\in\Cal O_L$. Fix a uniformizer $\pi\in F$ and let $K|F$ be the Galois extension of degree $p$ which we get by adjoining a $p$-th root $\psi$ of $\pi$. Let $K(E)$ be the field of definition of the $\overline K$-valued points of $E[p]$ where $E$ is an admissible elliptic curve over $F$. Let $e'=\text{ord}_K(p)$ be the absolute ramification index of $K$ and let $m'=e'/(p-1)$.
\enddefinition
\proclaim{Lemma 4.9} We have $F(E)^*\subseteq U_{K(E)}^{({m'p\over 2})}
(K(E)^*)^p$.
\endproclaim
\definition{Proof} Let $\theta:\bold k_{F(E)}\rightarrow\Cal O_{F(E)}$ be a multiplicative system of representatives whose existence is guaranteed by Proposition 8 on page 35 of [16], and let $x\in F(E)^*$ be arbitrary. Since $\pi$ is also a uniformizer in $F(E)$ we may write $x$ in the form:
$$x=\pi^{\text{ord}_{F(E)}(x)}\left(\theta(x_0)+\theta(x_1)\pi+\cdots+\theta(x_k)\pi^k+\cdots\right)$$
for some $x_0,x_1,\ldots\in\bold k_{F(E)}$. For every $k\in\Bbb N$ let $y_k\in\bold k_{F(E)}$ be the $p$-th root of $x_k$ and let
$$y=\psi^{\text{ord}_{F(E)}(x)}
\left(\theta(y_0)+\theta(y_1)\psi+\cdots+\theta(y_k)\psi^k+\cdots\right)
\in K(E)^*.$$
Clearly $xy^{-p}\in\Cal O^*_{K(E)}$. Since $p\in\pi^{{m'p\over2}}\Cal O_{K(E)}$ we have:
$$(a_0+a_1+\cdots+a_k+\cdots)^p\equiv
(a_0^p+a_1^p+\cdots+a_k^p+\cdots)\mod\pi^{{m'p\over2}}\Cal O_{K(E)}$$
for any convergent series $a_0,a_1,\ldots,a_k,\ldots\in\Cal O_{K(E)}$, and hence
$$\split xy^{-p}=&
\left(\theta(x_0)+\cdots+\theta(x_k)\pi^k+\cdots\right)
\cdot\left(
\theta(y_0)+\cdots+\theta(y_k)\psi^k+\cdots\right)^{-p}\\
\equiv&
\left(\theta(x_0)+\cdots+\theta(x_k)\pi^k+\cdots\right)\cdot
\left(\theta(y_0)^p+\cdots+\theta(y_k)^p\pi^k+\cdots\right)^{-1}\\
\equiv&\ 1\mod\pi^{{m'p\over2}}\Cal O_{K(E)},\endsplit$$
where we used that $\theta(y_k)^p=\theta(x_k)$ for every $k\in\Bbb N$. The claim is now clear.\qed
\enddefinition
\definition{Proof of Theorem 1.1} As explained in [14] we only need to show the following; let $E$ be an elliptic curve defined over $F$ and let $c$ be a cohomology class in $H^1(F,E)[p]$. Then there is a finite solvable Galois extension $K|F$ such that the image of $c$ under the restriction map $r_{K|F}:H^1(F,E)\rightarrow H^1(K^{per},E)$ is zero. It is also explained in [14] that we only need to prove this claim when the action of $\Gamma_F$ on $E[p]$ does not factor through a solvable quotient. Therefore we may assume that the characteristic of $F$ is zero, and moreover we may assume by Proposition 4.7 that $E$ has good reduction over $F$ and it is admissible, by taking a suitable finite solvable Galois extension of $F$, if it is necessary. The claim now follows from Lemma 4.9 and Theorem 3.6 applied to the field $K$ introduced in Notation 4.8.\qed
\enddefinition

\heading 5. Extensions whose ramification does not disappear after any solvable extension
\endheading

\definition{Notation 5.1} For every $\Bbb F_p[SL_2(\Bbb F_p)]$-module $M$ we say that a group $G$ is $M$-type if it is the extension of $SL_2(\Bbb F_p)$ by $M$ and the action of $SL_2(\Bbb F_p)$ on $M$ via conjugation in $G$ is the given one. Given two Galois extensions $K|F$ and $L|F$ let $KL|F$ denote their composite. We continue to denote by $F$ a local field of characteristic zero with a residue field of characteristic $p>0$.
\enddefinition
\proclaim{Lemma 5.2} Assume that $p\geq5$ and let $L|F$ be a finite Galois extension whose Galois group is $M$-type for an $\Bbb F_p[SL_2(\Bbb F_p)]$-module $M$. Also suppose that the image of $I_F$ in $\text{\rm Gal}(L|F)$ is an $\Bbb F_p[SL_2(\Bbb F_p)]$-submodule $N\leq M$ which is isomorphic to $\Bbb F_p^2$ with its usual $SL_2(\Bbb F_p)$-action. Then for every finite solvable Galois extension $P|F$ the image of $I_P$ in $\text{\rm Gal}(LP|P)\leq \text{\rm Gal}(L|F)$ is also $N$.
\endproclaim
\definition{Proof} Note that Gal$(LP|P)$ is a normal subgroup of Gal$(L|F)$ whose quotient subgroup is isomorphic to Gal$(P|F)$, and hence the image of Gal$(LP|P)$ in the quotient $SL_2(\Bbb F_p)$ of $\text{Gal}(L|F)$ by $M$ is also a normal subgroup with a solvable quotient. Since the only such subgroup of $SL_2(\Bbb F_p)$ is itself we get that Gal$(LP|P)$ surjects onto the quotient $SL_2(\Bbb F_p)$ of $\text{Gal}(L|F)$. This means that we may reduce to the case when Gal$(P|F)$ is a group whose order is a prime $l$ by induction on the order of Gal$(P|F)$. The quotient group $I_F/I_P$ is isomorphic to Gal$(P^{un}|F^{un})$ where $F^{un}$ and $P^{un}$ are the maximal unramified extensions of $F$ and $P$, respectively. Since Gal$(P^{un}|F^{un})$ is a subgroup of Gal$(P|F)$ we get that the index of $I_P$ in $I_F$ divides $l$. Therefore the image $R$ of $I_P$ in $\text{\rm Gal}(LP|P)$ is a normal subgroup of $N$ of index dividing $l$. Moreover $R$ is a normal subgroup in Gal$(LP|P)$ so it is invariant under the action of $SL_2(\Bbb F_p)$. Because the action of $SL_2(\Bbb F_p)$ on $N$ is irreducible we get that $R$ is equal to $N$.\qed
\enddefinition
In the next three lemmas we will assume that $F$ contains the $p$-th roots of unity. Let $m=e/(p-1)$ as above.
\proclaim{Lemma 5.3} Every extension $K|F$ which we get by adjoining a $p$-th root of any $\alpha\in U_F^{(mp)}$ is unramified.
\endproclaim
\definition{Proof} We may assume without the loss of generality that $\alpha$ is not a $p$-th root in $F$. Let $\beta\in K$ be a $p$-th root of $\alpha$ and fix a uniformizer $\pi\in F$. Then $\gamma=\pi^{-m}(\beta-1)$ generates the degree $p$ extension $K|F$. The minimal polynomial of this element is:
$$f(x)=\pi^{-mp}\left((\pi^mx+1)^p-\alpha\right)=\pi^{-mp}(1-\alpha)+\sum_{k=1}^{p-1}\left({p\atop k}\right)\pi^{m(k-p)}x^k+x^p\in F[x].$$
This is a monic polynomial whose constant term $\pi^{-mp}(1-\alpha)$ is in $\Cal O_F$ by assumption. For every integer $k$ with $1\leq k\leq p-1$ we have:
$$\text{ord}_F(\left({p\atop k}\right)\pi^{m(k-p)})=m(p-1)+m(k-p)= m(k-1)\geq0,$$
so $f(x)\in\Cal O_F[x]$. Let $L|F$ be the unique unramified Galois extension such that $\bold k_L$ is the splitting field of $f(x)$ modulo $\pi\Cal O_F$. Since $f'(x)\equiv p\pi^{m(1-p)}\not\equiv0\mod\pi\Cal O_F$ the polynomial $f(x)$ splits in $L$ by Hensel's lemma. Therefore $K=F(\gamma)$ is contained in $L$, and hence it is an unramified extension of $F$ as we claimed.\qed
\enddefinition
\proclaim{Lemma 5.4} We have $(F^*)^p\cap U_F^{(mp-1)}\subseteq U_F^{(mp)}$. 
\endproclaim
\definition{Proof} Fix a uniformizer $\pi\in\Cal O_F$ and let $x\in F^*$ such that $x^p\in U_F^{(mp-1)}$. Clearly we have $x\in U^{(1)}_F$ and hence there is a positive integer $n$ and an $x_0\in\Cal O^*_F$ such that $x=1+\pi^nx_0$. By the binomial theorem:
$$x^p=(1+\pi^nx_0)^p=1+p\pi^nx_0+\cdots+\left({p\atop k}\right)\pi^{nk}x_0^k+\cdots+\pi^{np}x_0^p.\tag5.4.1$$
Assume first that $n<m$. Then for every integer $k$ with $1\leq k\leq p-1$ we have:
$$\text{ord}_F(\pi^{np}x_0^p)=pn<(p-1)m+nk=
\text{ord}_F(\left({p\atop k}\right)\pi^{nk}x_0^k),$$
and hence $x^p\notin U_F^{(mp-1)}$. This is a contradiction, so $n\geq m$. In this case
$$pm\leq \text{ord}_F(\pi^{np}x_0^p)\text{ and }pm\leq
\text{ord}_F(\left({p\atop k}\right)\pi^{nk}x_0^k),$$
and hence $x^p\in U^{(mp)}_F$ as we claimed.\qed
\enddefinition
For every group $G$ and for every $\Bbb F_p[G]$-module $N$ let $N^{\vee}$ denote its dual $\Bbb F_p[G]$-module $\text{Hom}_{\Bbb F_p}(N,\Bbb F_p)$. (Recall that for a left $\Bbb F_p[G]$-module $N$ we define the left $G$-multiplication on $N^{\vee}$ by the formula $g\lambda(x)=\lambda(g^{-1}x)$ for every $g\in G,\lambda\in\text{Hom}_{\Bbb F_p}(N,\Bbb F_p)$ and $x\in N$.)
\proclaim{Lemma 5.5} Let $K|F$ be a finite unramified Galois extension with Galois group $SL_2(\Bbb F_p)$. Let $N\leq\bold k_K$ be an irreducible finite $\Bbb F_p[\text{\rm Gal}(\bold k_K|\bold k_F)]=\Bbb F_p[\text{\rm Gal}(K|F)]$-module. Then there is a finite Galois extension $L|F$ such that
\roster
\item"$(i)$" the group $\text{\rm Gal}(L|F)$ is $M$-type for a finite $\Bbb F_p[SL_2(\Bbb F_p)]$-module $M$ which contains $N^{\vee}$,
\item"$(ii)$" the field $L$ contains $K$ and $M\cong\text{\rm Gal}(L|K)$,
\item"$(iii)$" the image of $I_F$ in $M\leq\text{\rm Gal}(L|F)$ is $N^{\vee}$.
\endroster
\endproclaim
\definition{Proof} Fix a uniformizer $\pi\in\Cal O_F$. Then the homomorphism $\Cal O_K\rightarrow U_K^{(mp-1)}$ given by the rule $x\mapsto 1+x\pi^{mp-1}$ maps $\pi\Cal O_F$ onto $U_F^{(mp)}$ and hence induces an isomorphism:
$$\bold k_K\cong \Cal O_K/\pi\Cal O_K\cong U_K^{(mp-1)}
/U_K^{(mp)}.$$
This isomorphism is Gal$(K|F)$-equivariant where we consider $\bold k_K$ as a Gal$(K|F)$-module via the natural action of Gal$(\bold k_K|\bold k_F)$ and we equip $U_K^{(mp-1)}/U_K^{(mp)}$ with the quotient Gal$(K|F)$-module structure. By applying Lemma 5.4 to the local field $K$ we get that there is a natural quotient map:
$$q:U_K^{(mp-1)}/ U_K^{(mp-1)}\cap(K^*)^p
\rightarrow U_K^{(mp-1)}/ U_K^{(mp)}.$$
Choose a finite set $\overline N\subset U_K^{(mp-1)}/ U_K^{(mp-1)}\cap(K^*)^p$ such that the image of $\overline N$ under $q$ is $N$. Let $M_0\subseteq U_K^{(mp-1)}/U_K^{(mp-1)}\cap(K^*)^p$ be the $\Bbb F_p[\text{Gal}(K|F)]$-module generated by $\overline N$. Since $M_0$ is the $\Bbb F_p$-span of all the $\text{Gal}(K|F)$-conjugates of $\overline N$, and since that set is finite, we get that $M_0$ is a finite $\Bbb F_p[\text{Gal}(K|F)]$-module. Let $\overline M\subseteq U_K^{(mp-1)}$ be a set of representatives of $M_0$ and let $L|K$ be the finite Galois extension which we get by adjoining the $p$-th roots of $\overline M$. This extension is independent of the choice of $\overline M$ and since $M_0$ is invariant under the action of Gal$(K|F)$ the extension $L|F$ is also Galois. By Kummer theory Gal$(L|K)$ is isomorphic to $M=M_0^{\vee}$ as a Gal$(K|F)$-module. The image of $M_0$ with respect to $q$ is $N$ and hence $M_0$ has a quotient $\Bbb F_p[\text{Gal}(K|F)]$-module isomorphic to $N$. By duality $M$ has a $\Bbb F_p[\text{Gal}(K|F)]$-submodule isomorphic to $N^{\vee}$. So properties $(i)$ and $(ii)$ hold for $L|F$.

Let $M_0'$ be the kernel of the restriction of $q$ onto $M_0$ and let $\overline M'\subseteq U_K^{(mp)}$ be a set of representatives of $M_0'$. Let $L'|K$ be the the finite Galois extension which we get by adjoining the $p$-th roots of $\overline M'$. By Lemma 5.3 the extension $L'|K$ is unramified. Since the extension $K|F$ is also unramified we get that $L'|F$ is an unramified extension, too. Because $N^{\vee}\leq\text{Gal}(L|F)$ is equal to the Galois group Gal$(L|L')$, the image of $I_F$ in $\text{Gal}(L|F)$ is contained in $N^{\vee}$. By Lemma 5.4 there is not any unramified extension of $K$ where any element of $U^{(mp-1)}_K-U^{(mp)}_K$ could be a $p$-th root. Therefore the extension $L|K$ is ramified and so the image of $I_F$ in $N^{\vee}$ is non-trivial. Because $N^{\vee}$ is the dual of an irreducible $\Bbb F_p[\text{Gal}(K|F)]$-module, it is also irreducible and hence the image of $I_F$ must be $N^{\vee}$ itself.\qed
\enddefinition
\proclaim{Proposition 5.6} Let $\bold f$ be an infinite perfect field of characteristic $p$ whose Brauer group has trivial $2$-torsion and let $\bold k|\bold f$ be a finite Galois extension. Let $N$ be an $\Bbb F_p[\text{\rm Gal}(\bold k|\bold f)]$-module which has dimension two as a vector space over $\Bbb F_p$. Then there is a $\Bbb F_p[\text{\rm Gal}(\bold k|\bold f)]$-submodule of $\bold k$ which is isomorphic to $N$.
\endproclaim
\definition{Proof} Note that it will be enough to show that there is an $\Bbb F_p[\text{\rm Gal}(\overline{\bold f}|\bold f)]$-submodule $M$ of $\overline{\bold f}$ which is isomorphic to $N$ as an $\Bbb F_p[\text{\rm Gal}(\overline{\bold f}|\bold f)]$-module. In fact in this case the action of $\text{\rm Gal}(\overline{\bold f}|\bold f)$ on the extension of $\bold f$ generated by $M$ factors through Gal$(\bold k|\bold f)$. Therefore $M$ will be a submodule of $\bold k$ by the fundamental theorem of Galois theory. Let $A$ denote the polynomial ring $\Bbb F_p[t]$. We may consider $\bold f$ as an extension of the residue field $\Bbb F_p=A/(t+1)$ of the prime ideal $(t+1)\triangleleft A$ and hence we may talk about Drinfeld $A$-modules of rank $2$ of characteristic $(t+1)$ with coefficients in $\bold f$. For such a Drinfeld $A$-module $\phi$ its $t$-torsion group scheme $\phi[t]$ is \'etale and hence $\phi[t](\overline{\bold f})$ is a $\Bbb F_p[\text{\rm Gal}(\overline{\bold f}|\bold f)]$-submodule of $\overline{\bold f}$ which has dimension two as a vector space over $\Bbb F_p$.

Let $Y(t)$ be the Drinfeld modular curve over $\Bbb F_p=A/(t+1)$ which parameterizes Drinfeld $A$-modules of rank $2$ of characteristic $(t+1)$ with a full $t$-level structure. By fixing an $\Bbb F_p$-basis of $N$ we get a continuous Galois representation $\rho:\text{Gal}(\overline{\bold f}|\bold f)\rightarrow GL_2(\Bbb F_p)$. Let $Y_{\rho}$ be the twist of $Y(t)$ with respect to $\rho$ via the natural action of $GL_2(\Bbb F_p)$ on $Y(t)$. The $\bold f$-valued rational points of $Y_{\rho}$ correspond to Drinfeld $A$-modules of rank $2$ of characteristic $(t+1)$ with coefficients in $\bold f$ such that the $\text{Gal}(\overline{\bold f}|\bold f)$-module $\phi[t](\overline{\bold f})$ is isomorphic to $N$. Since the smooth compatification of the affine curve $Y(t)$ is geometrically irreducible of genus zero, the same holds for $Y_{\rho}$. As the $2$-torsion of the Brauer group of $\bold f$ is zero we get that $Y_{\rho}$ has a Zariski-dense set of $\bold f$-valued points.\qed
\enddefinition
\definition{Proof of Theorem 1.2} Let $p\geq5$ and let $\bold f$ be the perfection of a function field of transcendence degree one over an algebraically closed field of characteristic $p$. Then there is a finite Galois extension $\bold k|\bold f$ with Galois group $SL_2(\Bbb F_p)$ by Harbater's theorem (see Corollary 1.5 of [7] on page 284). Let $F$ be the field which we get by adjoining the $p$-th roots of unity to the fraction field of the ring of Witt vectors of $\bold f$ and let $K$ be the unique unramified extension of $F$ with residue field $\bold k$. By Tsen's theorem the Brauer group of $\bold f$ is trivial and hence there is a $\text{\rm Gal}(\bold k|\bold f)$-submodule $N\leq\bold k_K$ which is isomorphic to $\Bbb F_p^2$ with its usual $SL_2(\Bbb F_p)$-action by Proposition 5.6. Because $\Bbb F_p^2$ is self-dual as an $\Bbb F_p[SL_2(\Bbb F_p)]$-module we may apply Lemma 5.5 to the extension $K|F$ and the module $N$ to get a finite Galois extension $L|F$ which satisfies the conditions in Lemma 5.2. The claim now follows from this lemma.\qed
\enddefinition

\heading 6. Abelian varieties which do not become semi-stable after any solvable extension
\endheading

\definition{Definition 6.1} For every algebraic group $T$ over a field $F$ which is a split torus, let $C(T)$ denote its group of cocharacters. Then $C(T)$ is a free and finitely generated abelian group whose rank is equal to the dimension of $T$ over $F$. The group $T(F)$ of $F$-valued points of $T$ is canonically isomorphic to $F^*\otimes C(T)$. Now let $F$ be a field complete with respect to a discrete valuation $v:F^*\rightarrow\Bbb Z$. A subgroup $\Lambda$ of $T(F)=F^*\otimes C(T)$ is called a discrete lattice if the restriction of the homomorphism $v\otimes1:F^*\otimes C(T)\rightarrow
\Bbb Z\otimes C(T)$ to $\Lambda$ is injective and it has finite cokernel. In this case the quotient $T/\Lambda$ exists in the category of rigid analytic spaces and it is a proper rigid analytic group such that the quotient map $T\rightarrow T/\Lambda$ is a homomorphism of rigid analytic groups (see pages 324--325 of [4]). Let $\text{End}(T,\Gamma)$ denote the ring of endomorphisms of the algebraic group $T$ over $F$ such that $\phi(\Gamma)\subseteq\Gamma$. The operation of forming quotients induces an injective homomorphism
$$h_{T,\Lambda}:\text{End}(T,\Gamma)\longrightarrow\text{End}(T/\Gamma)$$
where the latter is the ring of rigid analytic endomorphisms from the rigid
analytic group $T/\Lambda$.
\enddefinition
\proclaim{Theorem 6.2} The homomorphism $h_{T,\Lambda}$ is an isomorphism.
\endproclaim
\definition{Proof} This is Satz 5 of [3] on page 33.\qed
\enddefinition
\proclaim{Theorem 6.3} Let $F$ be a local field, let $T$ be a split torus over $F$ and let $\Lambda\subset T(F)$ be a discrete lattice. Then the quotient $T/\Lambda$ is isomorphic to the rigid analytic variety underlying an abelian variety over $F$ if and only if there is a homomorphism $\lambda$ from $\Lambda$ to the character group $\text{\rm Hom}(T,\Bbb G_m)$ of $T$ such that the bilinear map:
$$(\alpha,\beta)\mapsto\lambda(\alpha)(\beta):\Lambda\times\Lambda\rightarrow F^*$$
is symmetric and $\text{\rm ord}_F((\alpha,\alpha))>0$ whenever $1\neq\alpha\in\Lambda$. 
\endproclaim
\definition{Proof} This is Theorem 5 of [4] on page 338.\qed
\enddefinition
\definition{Notation 6.4} For every field $F$ and for every commutative group scheme $B$ over $F$ let End$(B)$ denote the ring of endomorphisms of $B$ as a group scheme over $F$. Moreover let Aut$(B)\subseteq\text{End}(B)$ denote the group of automorphisms of $B$ over $F$. Let $F$ be a local field, let $T$ be a split torus over $F$ and let $\Lambda\subset T(F)$ be a discrete lattice such that the quotient $T/\Lambda$ isomorphic to the rigid analytic variety underlying an abelian variety $B$ over $F$. By GAGA (see Theorem 2.8 of [9] on page 349) the ring of endomorphisms of $B$ as a group scheme over $F$ and the ring of rigid analytic endomorphisms of the rigid analytic group $T/\Lambda$ are the same. Let $g_B:\text{End}(B)\rightarrow\text{End}(C(T))$ denote the ring homomorphism which is the composition of
$$h_{T,\Lambda}^{-1}:\text{End}(B)=\text{End}(T/\Lambda)
\longrightarrow\text{End}(T,\Lambda)$$
and the forgetful map:
$$\text{End}(T,\Lambda)\longrightarrow\text{End}(T)=\text{End}(C(T)).$$ 
Finally let $f_B$ denote the $\Bbb Q$-linear representation
$$(g_B)|_{\text{Aut}(B)}\otimes\text{id}_{\Bbb Q}:\text{Aut}(B)\rightarrow GL(C(T)\otimes\Bbb Q).$$
\enddefinition
\proclaim{Lemma 6.5} Let $F$ be a local field and let $\rho:G\rightarrow GL(V)$ be a representation of the finite group $G$ on a finite dimensional vector space $V$ over $\Bbb Q$. Then there is a split torus $T$ over $F$ and a discrete lattice $\Lambda\subset T(F)$ such that the following holds:
\roster
\item"$(i)$" the quotient $T/\Lambda$ isomorphic to the rigid analytic variety underlying an abelian variety $B$ over $F$,
\item"$(ii)$" there is an injective homomorphism $\sigma:G\rightarrow \text{\rm Aut}(B)$,
\item"$(iii)$" the representation $(f_B)\circ\sigma$ is isomorphic to $\rho$.
\endroster
\endproclaim
\definition{Proof} Note that there is a $G$-submodule $\Gamma<V$ which is a finitely generated, free $\Bbb Z$-module of rank dim$(V)$ and spans $V$ as a vector space over $\Bbb Q$. Let $T$ be a split torus over $F$ such that $C(T)$ is equal to $\Gamma$. Choose a uniformizer $\pi\in F^*$ and let $\Lambda<F^*\otimes C(T)$ be the abelian group generated by the set $\{\pi\otimes\gamma|\gamma\in C(T)\}$. Let $v:F^*\rightarrow\Bbb Z$ be the valuation on $F$ normalised such that $v(\pi)=1$. Under the homomorphism $v\otimes\text{id}_{C(T)}$ the group $\Lambda$ maps isomorphically onto its image which is $\Bbb Z\otimes C(T)$. Therefore $\Lambda$ is a discrete lattice in $T(F)$.

Let $\langle\cdot,\cdot\rangle:V\times V\rightarrow\Bbb Q$ be a symmetric positive definite bilinear form. We may assume, by multiplying by a positive integer, if it is necessary, that $\langle\cdot,\cdot\rangle$ takes integer values on $C(T)=\Gamma$. Let $\theta:C(T)\rightarrow\text{Hom}(C(T),\Bbb Z)$ be the homomorphism given by the rule $\theta(\alpha)(\beta)=\langle\beta,\alpha\rangle$ for every $\alpha,\beta\in C(T)$. Let $\lambda:\Lambda\rightarrow\text{Hom}(T,\Bbb G_m)$ denote the composition:
$$\CD\Lambda @>{v\otimes\text{id}_{C(T)}}>>\Bbb Z\otimes C(T)\cong C(T)@>\theta>>
\text{Hom}(C(T),\Bbb Z)\cong\text{Hom}(T,\Bbb G_m)\endCD.$$
The map $\lambda$ satisfies the conditions of Theorem 6.3, and hence the quotient $T/\Lambda$ is isomorphic to the rigid analytic variety underlying an abelian variety $B$ over $F$. Note that the natural action of $G$ on $F^*\otimes C(T)$ induced by the $G$-module structure on $C(T)=\Gamma$ leaves the subgroup $\Lambda$ invariant. Therefore the natural action of $G$ on $T$ descends down to the quotient $T/\Lambda$, that is, there is an injective homomorphism $\sigma:G\rightarrow \text{\rm Aut}(B)$ by GAGA. The objects $T,\Lambda,\sigma$ obviously satisfy condition $(iii)$ above.\qed
\enddefinition
For every group $G$, for every representation $\rho:G\rightarrow GL(V)$ on a vector space $V$ over $\Bbb Q$, and for every prime number $l$ let $\rho_l:G
\rightarrow GL(V\otimes_{\Bbb Q}\Bbb Q_l)$ denote the $\Bbb Q_l$-linear extension of $\rho$.
\proclaim{Proposition 6.6} Let $F$ be a local field and let $\rho:\text{\rm Gal}(\overline F|F)\rightarrow GL(V)$ be a representation on a finite dimensional vector space $V$ over $\Bbb Q$ which is continuous with respect to the discrete topology on $GL(V)$. Then there is an abelian variety $A$ defined over $F$ such that for every prime number $l$ different from the characteristic of $F$ the $\text{\rm Gal}(\overline F|F)$-representation $V_l(A)$ has a quotient isomorphic to $\rho_l$. 
\endproclaim
\definition{Proof} Let $G$ be the image of $\text{\rm Gal}(\overline F|F)$ with respect to $\rho$ and by slight abuse of notation let $\rho$ denote the representation furnished by the inclusion $G\rightarrow GL(V)$. Let $T$ be a split torus over $F$ and let $\Lambda\subset T(F)$ be a discrete lattice which satisfies the properties of Lemma 6.4 with respect to the group $G$ and the representation $\rho$. As above let $B$ be an abelian variety over $K$ which is isomorphic to $T/\Lambda$ as a rigid analytic variety and let $\sigma:G\rightarrow\text{\rm Aut}(B)$ be a homomorphism which satisfies conditions $(ii)$ and $(iii)$ of Lemma 6.5.

Now let $l$ be a prime number different from the characteristic of $F$. The quotient map $i:T\rightarrow B$ induces an injective homomorphism $i_l:T_l(T)\rightarrow T_l(B)$ where $T_l(T)$ is the Tate module of the torus $T$. The cokernel of $i_l$ is isomorphic to $C(T)\otimes_{\Bbb Z}\Bbb Z_l$ equipped with the trivial Galois-action. Therefore there is a short exact sequence:
$$\CD0@>>>C(T)\otimes_{\Bbb Z}\Bbb Q_l(1)@>>> V_l(B)@>>>
C(T)\otimes_{\Bbb Z}\Bbb Q_l@>>>0\endCD\tag6.6.1$$
of $\text{\rm Gal}(\overline F|F)$-representations, where for every $\Bbb Q_l[\text{Gal}(\overline F|F)]$-module $W$ we let $W(1)$ denote the Tate twist of $W$, and we let $\text{Gal}(\overline F|F)$ act on $C(T)$ trivially. The composition of $\rho:\text{Gal}(\overline F|F)\rightarrow G$ and $\sigma$ furnishes a cohomology class in $H^1(F,\text{Aut}(B_{\overline F}))$. Let $A$ denote the twist of $B$ with respect to this class. Because the action of $G$ on $B$ lifts to $T$, its action on $V_l(B)$ respects the filtration of (6.6.1). Therefore there is a short exact sequence:
$$\CD0@>>>V\otimes_{\Bbb Q}\Bbb Q_l(1)@>>> V_l(A)@>>>
V\otimes_{\Bbb Q}\Bbb Q_l@>>>0\endCD\tag6.6.2$$
where Gal$(\overline F|F)$ acts on $V$ via $\rho$. The claim is now clear.\qed
\enddefinition
\proclaim{Proposition 6.7} There is a local field $F$ of characteristic zero with a perfect residue field and there is an abelian variety over $F$ which does not acquire semi-stable reduction over any solvable extension of $F$.
\endproclaim
\definition{Proof} By Theorem 1.2 there is a local field $F$ of characteristic zero with a perfect residue field and there is a finite Galois extension $K|F$ such that there is no solvable extension $L|F$ such that the extension $KL|K$ is unramified, where $KL$ is the composite of $K$ and $L$. Choose a faithful finite-dimensional $\Bbb Q$-linear representation $\rho:\text{Gal}(K|F)\rightarrow GL(V)$ and by slight abuse of notation let $\rho$ also denote the composition of the quotient map Gal$(\overline F|F)\rightarrow\text{Gal}(K|F)$ and $\rho$. Choose a prime number $l$ different from the characteristic of $\bold k_F$. By Proposition 6.6 there is an abelian variety $A$ over $F$ such that $V_l(A)$ has a quotient $W$ isomorphic to $\rho_l$ as a Gal$(\overline F|F)$-representation. 

We claim that $A$ does not acquire semi-stable reduction over any solvable extension of $F$. Assume the contrary and let $L|F$ a finite solvable Galois extension such that $A_L$ has semi-stable reduction. By Theorem 4.3 the action of $I_L$ on $V_l(A)$ is unipotent. Therefore the the action of $I_L$ on $W$ is unipotent, too. In particular the image of $I_L$ under $\rho_l$ is a pro-$l$ group, and hence the extension Gal$(KL|L)$ is tamely ramified. Therefore there is a finite solvable Galois extension $L'|L$ such that Gal$(KL'|L')$ is unramified by Lemma 4.6. The extension $L'|F$ is solvable as it is a tower of solvable extensions. This is a contradiction.\qed
\enddefinition
\definition{Proof of Corollary 1.3} By the proposition above there is a field $F$ of characteristic zero complete with respect to a discrete valuation whose residue field is perfect and there is an abelian variety $A$ over $F$ which does not acquire semi-stable reduction over any solvable extension of $F$. For any geometrically irreducible, smooth, projective curve $D$ we let Jac$(D)$ denote the Jacobian of $D$. By Faltings's trick (see Theorem 10.1 of [12] on page 198) there is a smooth, projective, geometrically irreducible curve $C$ over $F$ such that there is a surjective homomorphism $\phi:\text{Jac}(C)\rightarrow A$ of abelian varieties over $F$. We claim that $C$ does not acquire semi-stable reduction over any solvable extension of $F$. Assume the contrary and let $K|F$ a finite solvable Galois extension such that $C_K$ has semi-stable reduction. Then $\text{Jac}(C)_K=\text{Jac}(C_K)$ also has semi-stable reduction by Theorem 2.4 of [2] on page 89. Therefore the action of $I_K$ on $V_l(\text{Jac}(C)_K)$ is unipotent by Theorem 4.3 where $l$ is a prime number different from the characteristic of $\bold k_F$. Because the homomorphism $\phi_l:V_l(\text{Jac}(C)_K)\rightarrow V_l(A_K)$ induced by $\phi$ is surjective we get that the action of $I_K$ on $V_l(A_K)$ is unipotent, too. Therefore $A_K$ has semi-stable reduction by Theorem 4.3. This is a contradiction.\qed
\enddefinition


\Refs
\ref\no 1\by M.~\c{C}iperiani and A.~Wiles\paper Solvable points on genus one curves\jour Duke Math. J.\vol 142 \yr 2008\pages 381--464\endref

\ref\no 2\by P. Deligne and D. Mumford\paper The irreducibility of the space of curves of given genus\jour Publ. Math. Inst. Hautes \'Etudes Sci.\vol 36\yr 1969\pages 75--109\endref

\ref\no 3\by L. Gerritzen\paper \"Uber Endomorphismen nichtarchimedischer holomorpher Tori\jour Invent. Math.\vol 11\yr 1970\pages 27--36\endref

\ref\no 4\by L. Gerritzen\paper On non-archimedean representations
of abelian varieties\jour Math. Ann.\vol 196\yr 1972\pages
323--346\endref

\ref\no 5\by D. Gorenstein\book Finite groups\bookinfo second edition\publ Americal Mathematical Society\publaddr Providence, Rhode Island\yr 2007\endref

\ref\no 6\by A. Grothendieck and M. Raynaud\paper Mod\`eles de N\'eron et monodromie\inbook Groupes de monodromie en G\'eometrie Algebrique, I, II\publ Springer\bookinfo Lecture Notes in Math.\vol 288\yr 1972\pages 313--523\endref

\ref\no 7\by D. Harbater\paper Mock covers and Galois extensions\jour J. Algebra\vol 91\yr 1984\pages 281--293\endref

\ref\no 8\by Gy. K\'arolyi and A. P\'al\paper The cyclomatic number of connected graphs without solvable orbits\jour  J. Ramanujan Math. Soc.\paperinfo to appear\yr 2012\endref

\ref\no 9\by W. L\"utkebohmert\paper Formal-algebraic and rigid-analytic geometry\jour Math. Ann.\vol 286\yr 1990\pages 341--371\endref

\ref\no 10\by B. Mazur and L. Roberts\paper Local Euler characteristics\jour Invent. Math.\vol 9\yr 1970\pages 201--223\endref

\ref\no 11\by J. Milne\book \'Etale cohomology\publ Princeton University Press\publaddr Princeton, New Jersey\yr 1980\endref

\ref\no 12\by J. Milne\paper Jacobian varieties\inbook Arithmetic geometry (Storrs, Conn., 1984)\pages 167–-212\publ Springer\publaddr New York\yr 1986\endref

\ref\no 13\by F. Oort and J. Tate\paper Group schemes of prime order\jour Annales scientifiques de l'\'E.N.S, 4e s\'erie, t.3\yr 1970\pages 1--21\endref

\ref\no 14\by A. P\'al\paper Solvable points on projective algebraic curves\jour Canad. J. Math.\vol 56\yr 2004\pages 612--637\endref

\ref\no 15\by L. Roberts\paper The flat cohomology of group schemes of rank $p$\jour American Journal of Mathematics\vol 95\yr 1973\pages 688--702\endref

\ref\no 16\by J.--P. Serre\book Local fields\bookinfo GTM 67 (second corrected pringting, 1995)\publ Springer--Verlag\publaddr New York--Berlin\yr 1979\endref

\ref\no 17\by J. Tate\paper Finite flat group schemes\inbook Modular forms and Fermat’s last theorem (Boston, MA, 1995)\publ Springer-Verlag\publaddr Berlin-Heidelberg-New York\yr 1997\pages 121--154\endref
\endRefs
\enddocument